\documentclass[12pt]{article}
\usepackage[dvips]{color}
 \usepackage{amssymb}
 \usepackage{amsmath}
 \usepackage{amsfonts}
 \usepackage{mathrsfs}
\usepackage{amsmath, amsthm, amssymb}
 \newcommand{\be}{\begin{equation}}
 \newcommand{\ee}{\end{equation}}
 \newcommand{\bea}{\begin{eqnarray}}
 \newcommand{\eea}{\end{eqnarray}}
 \newcommand{\nn}{\nonumber}
 \newcommand{\rd}{\partial}

 \newtheorem{thm}{Theorem}
 \newtheorem{definition}{Definition}
 
 \newtheorem{lem}{Lemma}

 \newtheorem{Rema}{Remark}
 \makeatletter
 \def\blfootnote{\xdef\@thefnmark{}\@footnotetext}
 \makeatother

\begin{document}

 \begin{center}
    \font\titlerm=cmr10 scaled\magstep3
    \font\titlei=cmmi10 scaled\magstep4
    \font\titleis=cmmi7 scaled\magstep4
     \centerline{\titlerm  Analytic solution for one dimensional inverse heat}
      \centerline{\titlerm  }
     \centerline{\titlerm  conduction problem of semi-infinite bar}
     \vspace{1.5cm}

    \noindent{{\large
       \, Adel
       Kassaian\,\, and\, \, A. Haghany}}
  \end{center}

\begin{center}
    \font\titleis=cmmi7 scaled\magstep3
 \noindent{{ E-mails : a.kassaian@gmail.com } }
  \noindent{{ \& \,\, aghagh@cc.iut.ac.ir} }
  \end{center}

\vspace{0.5cm}
\begin{abstract}
We present analytical formula along with its existence theorem for
solution of inverse heat conduction problem of semi-infinite bar,
equivalent to a Volterra integral equation of first kind, as an
infinite series of fractional derivatives. The mathematical method
is based on some properties of function space $M[0,T]$  (proved
here) with respect to fractional integration and derivatives.
\end{abstract}
\vspace{0.5cm}


\vspace{0.5cm}

 \vskip 1em
%----------------------------------------------------------------------------------------------------------------

\section {Introduction}
In this paper we consider the problem of finding analytic solution
for one dimensional inverse heat conduction problem of {\em
semi-infinite} bar (Cauchy problem (\ref{shahla})). This problem can
turn into a problem of Volterra integral equation of first kind
(VIE1),

 \be \int^{t}_{0} \,K(x_0-x,t-\tau) \,u(x,\tau)\,d\tau
=g(t),\hspace{1cm} 0 \le x < x_0 ,\,\,\,0 < t \le T,\nn\ee
 where, \be K(x,t)= \begin{cases}\displaystyle{{\,\,x\over \Large{ {2\sqrt{\pi}\,
 t^{3\over2}}}}}\,\text{e}^{-{\,\,x^2\over \Large{ {4\, t}}}} & \mbox { if }\,\,\, 0 < t \le T, \\
0  & \mbox { if }\,\,\,  t = 0.
\end{cases}\label{ahmad}\ee

The concern of this paper is only on finding the analytical solution
along with an existence theorem for this problem. For regularization
of this specific ill-posed problem the reader is referred to
\cite{lam}. Although considerable amount of work has been done for
analytical solution, existence condition and regularization of one
dimensional inverse heat conduction problem of {\em finite} bar
(Cauchy problem (\ref{goli})), e.g. see \cite{hao}, the analytical
solution of Cauchy problem (\ref{shahla}) or equivalently the above
VIE1 problem, seems to be not treated in the
literature.\\

The method we use to solve this problem is based on properties of
function space $M[0,T]$; see definition \ref{defi} below. This
function space is proved to be a useful tool in analyzing linear
integro differential causal problems \cite{ij}. In this paper we
prove and use different properties of $M[0,T]$,
 regarding fractional derivatives and fractional integration. The importance of this function space
 to this
problems is due to the fact that $K({ x},t)$ (and all partial
derivatives $\rd^{i}_x K({ x},t)$) for every $T>0$ as a function of
$t$ belongs to $M[0,T]$ for all $x\ne0$. From this it follows that
$g(t)$ belongs to $M[0,T]$, when the solution $u(x,t)$ is assumed to
be a continuous function of $t$. In section \ref{hamd} of this paper
we prove the fact that $M[0,T]$ is closed under
 fractional integrations and derivatives. In particular we show in
part 1 of Theorem \ref{lmf4} that the Riemann-Liouville fractional
derivatives
 and the Caputo fractional derivatives to all order exist and they coincide
  (for every fractional order) as
  elements of $M[0,T]$. The relation (\ref{hop}) in part 2 of Theorem \ref{lmf4} is used in section \ref{heat}
  to solve the above VIE1
  which is a convolution type VIE1 with its kernel being an element of
  $M[0,T]$.\\

In the part \ref{samme3}, by using the results of section \ref{hamd}
we are able to deduce from the above VIE1 in Theorem \ref{lemf88}
that,

\be \rd_{x}\,u(x,t)|_{x=x_{0}}=-\rd^{1\over2}_{t}g(t),\nn \ee

where $\rd^{1\over2}_{t}g(t)=\,{}_{RL}D^{1\over
2}_{0\,t}\,g(t)=\,{}_{C}D^{1\over 2}_{0\,t}\,g(t)$. Using the result
above, the problem of semi-infinite bar can turn into the problem of
finite bar (equation (\ref{goli})) and by which we find the answer
as, \be u(x,t) = \sum_{n=0}^{\infty}\rd^{n\over 2}_{t}(g (t))
{{(x_0-x)^{n}}\over(n)!}. \nn\ee  In the part \ref{ans} we discuss
the convergence of this solution for $u(x,t)$ and derive an
existence theorem for this problem using Holmgren classes. Finally
in \ref{example}, a convergent example is provided for the solution
$u(x,t)$, through a special case of initial value problem of one
dimensional direct heat conduction equation of an {\em infinite bar}
($x\in(-\infty,\infty$)).

\section {Function space $M[0,T]$ and fractional derivatives}
\label{hamd}

For the rest of this paper we adopt the following notations and definitions.\\

 $\mathbb{N} :=\{1, 2, \ldots \}$ is the set of natural numbers.
 $\mathbb{Z}_{+} :=\{0,1, 2, \ldots \}$ is the set of non-negative
 integers.  ${C}^{n}[a,b]$ is the space of $n$
times continuously differentiable functions on $[a,b]$.
${C}^{\infty}[a,b]$ is the space of
smooth functions on $[a,b]$.\\

Let us recall the definition of function space $M[0,T]$ from
\cite{ij}.
\begin{definition}\label{defi}{ Vector space $M[0,T]$ is the space
of all functions $\phi(x)\in C^{\infty}[0,T]$ for which,
\be\phi(0)=0\hspace{1cm}{\text{and}}\hspace{1cm}
\frac{d^n}{dt^n}\phi(t)|_{t=0}=0,\hspace{0.5cm}\forall n\in
\mathbb{N}.\label{moon}\ee}
\end{definition}
We also say $\phi(t)\in M[0,+\infty]$ if $\phi(t)\in M[0,T]$ for
every $T>0$. One can show for example that the function $K({x},t)$
given by (\ref{ahmad}) as a function of $t$ is in $M[0,+\infty]$ for
 ${{x}}\ne{0}$. It can be also easily checked
that function space $M[0,T]$ satisfies the following conditions. \be
\phi(t)\in M[0,T]\,\,\Longrightarrow\,\, \rd_{t}\phi(t)\in
M[0,T]\,\, ,\hspace{0.5cm} \int^{t}_{0} d\tau \phi(\tau)\in
M[0,T].\label{hossein}\ee

\be \phi(t)\in M[0,T],\,\,\, P(t)\in C[0,T]\,\Longrightarrow\,\,
 \int^{t}_{0} d\tau\,
\,\phi(t-\tau) P(\tau)\in M[0,T].\label{hasan}\ee

From (\ref{hossein}) one finds, the function space $M[0,T]$ is
closed under actions of derivative, $\rd_{t}:\phi(t) \mapsto
\rd_{t}\phi(t)$, and integral, $\int^{t}_{0} d\tau:\phi(t) \mapsto
\int^{t}_{0} d\tau \phi(\tau)$ operators. In this paper, in many
steps we use this property of $M[0,T]$ in the integration by parts
without mentioning. Also by (\ref{hasan}) each $\phi(t)\in M[0,T]$
can be viewed as a map from $C[0,T]$ to $M[0,T]$ via Volterra
convolution operator.

Considering fractional integral operator defined as, \be
J^{\alpha}_{c,\,t}(P(t))={1\over \Gamma(\alpha)}\int^{t}_{c}d\tau
(t-\tau)^{\alpha-1}P(\tau),\hspace{0.75cm}
0<\alpha<1,\hspace{0.25cm} P(t)\in C[c,T], \nn\ee where $0\le c<t$,
in the following lemma we prove $M[0,T]$ is closed under action of
fractional integral operator $J^{\alpha}_{0,t}$ and also, the
operators $J^{\alpha}_{0,t}$ and $\rd_t$ commute on elements of
$M[0,T]$.

\begin{lem}\label{lemf77}{
If $\phi(t)\in M[0,T]$ and $P(t)\in C[0,T]$ then for any
$0<\alpha<1$ we have,
\bea J^{\alpha}_{0,\,t}(\phi(t))&\in& M[0,T],\nn\\
\rd_{t}\big(J_{0,t}^{\alpha}(\phi(t))\big)&=&J_{0,\,t}^{\alpha}\big(\rd_{t}\phi(t)
\big),\nn\\J^{\alpha}_{0,t}\big(\int^{t}_{0} \phi(t-\tau)
P(\tau)\,d\tau\big)&=&\int^{t}_{0}
J^{\alpha}_{0,t-\tau}(\phi(t-\tau))P(\tau) d\tau.
\label{rezvan}\eea}
\end{lem}

{\bf Proof}. For the first part we note, \be
J^{\alpha}_{0,\,t}(\phi(t))={1\over \Gamma(\alpha)}\int^{t}_{0}d\tau
(t-\tau)^{\alpha-1}\phi(\tau)= {1\over
\Gamma(\alpha)\alpha}\int^{t}_{0}d\tau
(t-\tau)^{\alpha}\rd_{\tau}\phi(\tau).\nn\ee

Thus $J^{\alpha}_{0,\,t}(\phi(t))|_{t=0}=0$ and
$J^{\alpha}_{0,\,t}(\phi(t))$ is continuous on $[0,T]$. Furthermore
by last line we have, \be \rd_{t}(\,J^{\alpha}_{0,\,t}(\phi(t))\,)=
{1\over \Gamma(\alpha)}\int^{t}_{0}d\tau
(t-\tau)^{\alpha-1}\rd_{\tau}\phi(\tau)={1\over
\Gamma(\alpha)\alpha}\int^{t}_{0}d\tau
(t-\tau)^{\alpha}\rd^2_{\tau}\phi(\tau).\nn\ee

So,\, $\rd_{t}(\,J^{\alpha}_{0,\,t}(\phi(t))\,)|_{t=0}=0$ and
$\rd_{t}(\,J^{\alpha}_{0,\,t}(\phi(t))\,)$ is continuous on $[0,T]$.
Continuing in the same way one can show, \be
\rd^{\,n}_{t}(\,J^{\alpha}_{0,\,t}(\phi(t))\,)={1\over
\Gamma(\alpha)\alpha}\int^{t}_{0}d\tau
(t-\tau)^{\alpha}\rd^{\,n+1}_{\tau}\phi(\tau).\label{khaleghi} \ee

Therefore,

\be\forall n\in \mathbb{Z}_{+}, \hspace{0.5cm}
\rd^{n}_{t}(\,J^{\alpha}_{0,\,t}(\phi(t))\,)|_{t=0}=0 \hspace{0.5cm}
\text{and} \hspace{0.5cm} J^{\alpha}_{0,\,t}(\phi(t)) \in
C^{\,n}[0,T].\nn\ee Thus,
$J^{\alpha}_{0,\,t}(\phi(t))\in M[0,T]$.\\

For the second part using (\ref{khaleghi}) we have, \bea
\rd_{t}(\,J^{\alpha}_{0,\,t}(\phi(t))\,)={1\over
\Gamma(\alpha)\alpha}\int^{t}_{0}d\tau
(t-\tau)^{\alpha}\rd^2_{\tau}\phi(\tau) &=& {1\over
\Gamma(\alpha)}\int^{t}_{0}d\tau
(t-\tau)^{\alpha-1}\rd_{\tau}\phi(\tau)\nn \\&=&
\,J^{\alpha}_{0,\,t}((\rd_{t}\phi(t))\,).\nn\eea

For the third part by taking $F(t)=(\int^{t}_{0} \phi(t-\tau)
P(\tau)\,d\tau\big)$, we have \be \rd_t F(t)=(\int^{t}_{0} \rd_t
\phi(t-\tau) P(\tau)\,d\tau\big)\nn\ee and also from (\ref{hasan}),
$F(t)\in M[0,T]$. Now by starting from the left hand side of
(\ref{rezvan}) we have,

\bea J^{\alpha}_{0,\,t}(F(t))&=& {1\over
\Gamma(\alpha)\alpha}\int^{t}_{0}d\tau
(t-\tau)^{\alpha}\rd_{\tau}F(\tau)\nn\\&=&{1\over
\Gamma(\alpha)\alpha}\int^{t}_{0}d\tau
(t-\tau)^{\alpha}\int^{\tau}_{0} \rd_\tau \phi(\tau-z) P(z)\,dz \nn
\\&=&{1\over
\Gamma(\alpha)\alpha}\int^{t}_{0}\,dz\, G(t,z)\,P(z)\,dz.
\hspace{1.5cm} (\,\ddagger\,) \nn\eea

In the last line we use Volterra repeated integral formula (that is
if $B(t,\tau)$ and $L(t,\tau)$ are continuous on
$\triangle=\{(t,\tau)\,|\, 0 \leq \tau\leq t\leq T \}$ and $Q(t)$ is
a continuous function on $[0,T]$ then $ \int^{t}_{0}d\tau
B(t,\tau)\int^{\tau}_{0}\,dz\, L(\tau,z)Q(z)=\int^{t}_{0}\,dz\,
H(t,z) Q(z)$, where $H(t,z)=\int^{t}_{z}\,d\tau\,
B(t,\tau)L(\tau,z))$,\, so $G(t,z)$ is given by,

\bea G(t,z)&=&\int^{t}_{z}\,d\tau\, (t-\tau)^{\alpha}\rd_\tau
\phi(\tau-z)={\alpha}\int^{t}_{z}d\tau (t-\tau)^{\alpha-1}
\phi(\tau-z) \nn\\ &=&{\alpha}\int^{t-z}_{0}d\tau'
(t-\tau'-z)^{\alpha-1} \phi(\tau')={\alpha}\Gamma(\alpha)
\,\,J_{0,\,t-z}^{\,\alpha } (\phi(t-z)).\nn \eea

Thus by insetting the expression for $G(x,z)$ into $(\,\ddagger\,)$
we have. \be J^{\alpha}_{0,\,t}(F(t))=\int^{t}_{0}\,dz\,
J_{0,\,t-z}^{\,\alpha } (\phi(t-z))\,P(z)\,dz.\nn \ee  \qed \\

For any $\nu> 0$, with $[\nu]=m$ (that is, $m\le\nu<m+1$ and
$m\in\mathbb{Z}_{+}$), the Riemann-Liouville fractional derivative
and the Caputo fractional derivative are defined for $P(t)\in
C^{\infty}[c,T]$ as (e.g. see \cite{pold}), \bea
{}_{RL}D^{\nu}_{c,\,t}
P(t)=\rd_{t}^{m+1}\big({J}_{c,\,t}^{m-\nu+1}(P(t))\big),\label{luvil}\\
{}_{C}D^{\nu}_{c,\,t}
P(t)=(J_{c,\,t}^{m-\nu+1}(\rd_{t}^{m+1}P(t)).\label{caputo}\eea

Now we state the following theorem regrading fractional derivatives
${}_{RL}D^{\nu}_{0,\,t}$ and ${}_{C}D^{\nu}_{0,\,t}$ of elements of
$M[0,T]$.

\begin{thm}\label{lmf4}{1-If $\phi(t)\in M[0,T]$ then for any $\nu> 0$, with $[\nu]=m$, the Riemann-Liouville
fractional derivative
 and the Caputo fractional derivative exist and coincide in $M[0,T]$, in
 other word,

\bea \rd^{\,\nu}_{0,\,t}\,\phi(t)  &&=: \,\, {}_{RL}D^{\nu}_{0,\,t}
\phi(t)={}_{C}D^{\nu}_{0,\,t} \phi(t)\nn\\&&= \,\,{{\rd_{t}^{m+1}
\int^{t}_{0}\,ds\, (t-s)^{m-\nu}\phi(s)}\over \Gamma(m-\nu+1)}\nn\\
&&={{\int^{t}_{0}\,ds\, (t-s)^{m-\nu}\rd_{s}^{m+1}\phi(s)}\over
\Gamma(m-\nu+1)}\in M[0,T].\label{hope}\eea 2-For any $P(x)\in
C[0,T]$ and $\phi(x)\in M[0,T],$  \be
\rd^{\,\nu}_{0,\,t}\int^{t}_{0}
\phi(t-\tau)P(\tau)d\tau=\int^{t}_{0}
\big(\rd^{\,\nu}_{0,\,t-\tau}\phi(t-\tau)\big)P(\tau)d\tau.\label{hop}\ee}
\end{thm}

{\bf Proof.} Part 1 is the consequence of the fact that $M[0,T]$ is
closed under actions of both derivative $\rd_{t}$ and fractional
integration $J^{\alpha}_{0,\,t}$\, and also the fact that derivative
operator $\rd_{t}$ and fractional integration operator
$J^{\alpha}_{0,\,t}$ commute.\\

For part 2 we have,

\bea \rd^{\,\nu}_{0,\,t}\int^{t}_{0}
\phi(t-\tau)P(\tau)d\tau=\rd^{m+1}_{t}\Big(J^{m-\nu+1}_{0,\,t}\big(\int^{t}_{0}
\phi(t-\tau) P(\tau)\,d\tau\big)\Big)\nn\\=
\rd^{m+1}_{t}\int^{t}_{0}
J^{{m-\nu+1}}_{0,\,t-\tau}(\phi(t-\tau))P(\tau) d\tau.
\nn\hspace{1cm}(\diamond)\eea

By Lemma \ref{lemf77} we have $J^{{m-\nu+1}}_{0,\,t}(\phi(t))\in
M[0,T]$, it follows that, \be \forall n\in \mathbb{Z}_{+},
\hspace{1cm} \rd_{\,t}^{\,n}
\big(J^{{m-\nu+1}}_{0,\,t-\tau}\phi(t-\tau)\big)|_{t=\tau}=0, \nn\ee
therefore in ($\diamond$), $\rd^{m+1}_{t}$ can go inside the
integral, \bea \rd^{\,\nu}_{0,\,t}\int^{t}_{0}
\phi(t-\tau)P(\tau)d\tau &=& \int^{t}_{0}
\rd^{m+1}_{t}\big(J^{{m-\nu+1}}_{0,\,t-\tau}(\phi(t-\tau))\big)P(\tau)
d\tau \nn\\&=&\int^{t}_{0}
\big(\rd^{\,\nu}_{0,\,t-\tau}\phi(t-\tau)\big)P(\tau)d\tau.\nn\eea\qed

We state the following Lemma regarding some required properties of
$K(x,t)$ as our last needed mathematical tool for the following
section.
\begin{lem}\label{lemf6}
For $ K(x,t)$ given by (\ref{ahmad}) we have, \bea
 \rd_{x} K(x,t)&=&
\rd^{1\over2}_{0,\,t}K(x,t),\hspace{1.3cm}x>0,\,\, t\ge
0,\label{hop2}\\\lim_{x \to 0^+} \int^{t}_{0} d\tau K(x,t-\tau)
P(\tau)&=& P(t),\hspace{2.6cm} P(t)\in C[0,T].\label{hop3}\eea
\end{lem}

{\bf Proof.} Equations (\ref{hop2}) can be easily verified directly.
By applying $\rd^{1\over2}_{0,\,t}$ (or $\rd_{x}$) again on
(\ref{hop2}) one can reach to more familiar relation, $\rd_{t}
K(x,t)=\rd^2_{x} K(x,t)$. For equation (\ref{hop3}) for example see
\cite{JRC} lemma (4.2.3).\qed

\section {One dimensional inverse heat conduction problem of semi-infinite bar}\label{heat}

 In this section for convenience we denote,
 \be \,\rd^{\,\nu}_{0,\,t}: =\rd^{\,\nu}_{\,t},\hspace{0.5cm} \text{for}\hspace{0.2cm} \nu>0,\ee for elements of
$M[0,T]$.  Here we consider finding analytic solution for Cauchy
problem of determining a function $u(x,t)$ that satisfies, \bea
 \rd^2_x u(x,t) &=&\rd_t u(x,t), \hspace{1cm} 0\le x < \infty,\hspace{0.6cm} 0\le
 t \le T,
\nn \\ u(x_0,t) &=& g(t), \hspace{1.7cm} x_0 \in (0,\infty),
\hspace{0.5cm} 0\le
 t \le T,\label{shahla}\\ u(x,0) &=& 0, \hspace{2.1cm}
0\le x<\infty.\nn\eea

The equation (\ref{shahla}) can be considered as the formulation of
one dimensional inverse heat problem for retrieving $u(x,t)$, the
temperature at time $t$ ($0\le
 t \le T$) and at point $x$ ($x \ge 0$)
 along {\em semi-infinite} conducting bar ($x\in[0,\infty)$) on
non-negative $x$-axis which is initially at constant temperature
zero ($u(x,0) = 0$ for $x \in [0,\infty)$). The bar is insulated all
the way except at $x=0$ (where it is subjected to unknown heat
source) and our given data, \, $g(t)=u(x_0,t)$, is the temperature
function at specific point $x_{0}\in(0,\infty)$ for the time period,
$\,\,0\le
 t \le T$.\\

For $u(x,t)$, the solution of this problem, it can be proved  (e.g.
see \cite{JRC}),

\be \int^{t}_{0} \,K(x'-x,t-\tau) \,u(x,\tau)\,d\tau
=u(x',t),\hspace{1cm}x  < x'\nn,\ee

where $K(x,t)$ is given by (\ref{ahmad}). This leads to the
following equivalent problem of Volterra integral equation of first
kind  to equation (\ref{shahla}), for $0 < x \le x_0 $,

 \be \int^{t}_{0} \,K(x_0-x,t-\tau) \,u(x,\tau)\,d\tau
=g(t),\hspace{1cm} 0 \le x < x_0\,\,\,0\le
 t \le T \label{babol},\ee

and the following answer for $x_0 < x$,

\be u(x,t)=\int^{t}_{0} \,K(x-x_0,t-\tau)
\,g(\tau)\,d\tau,\hspace{1cm}x_0 < x,\,\,\,0\le
 t \le T\nn\ee

The function $K(x,t)$  (for $x\ne0$) as a function of $t$ is smooth
on $[0,T]$ and satisfies,

\be \forall n\in \mathbb{Z}_{+}, \hspace{1cm} \rd_{\,t}^{\,n}
K(x,t)|_{t=0}=0, \hspace{1cm} x\in \mathbb{R}-\{0\}. \nn\ee

where $\mathbb{Z}_{+}=\{0,1,2,...\}$. In equation (\ref{babol}) if
we assume that $u(x,t)$, as a function of $t$, is a continuous
function on $[0,T]$, then it follows (by differentiating
(\ref{babol})) that $g(t)$ is smooth on $[0,T]$ and it also
satisfies, \be \forall n\in \mathbb{Z}_{+}, \hspace{1cm}
\rd_{\,t}^{\,n} g(t)|_{t=0}=0, \nn\ee

thus $g(t)\in M[0,T]$.

\subsection{Analytic Solution}
\label{samme3}

Let us first begin by stating the problem and solution of analogous
case of inverse heat problem for {\em finite bar} ($x\in[0,l]$),

 \bea
 \rd^2_x w(x,t) &=&\rd_t w(x,t), \hspace{1cm} 0 \le x \le l,\hspace{0.5cm} 0\le
 t \le T,
\nn \\ w(l,t) &=& h(t), \hspace{1.8cm} 0\le
 t \le T, \label{goli}\\ \rd_x
w(x,t)|_{x=l} &=& f(t), \hspace{1.8cm} 0\le
 t \le T. \nn\eea Equation
(\ref{goli}) is the formulation of inverse heat conduction problem
of retrieving $w(x,t)$, the temperature at point $x$ ($0 \le x \le
l$) and time $t$ $(0\le
 t \le T)$ of a finite bar (located between points $x=0$
and $x=l$ of $x$-axis, insulated for $0< x < l$ and is subjected to
unknown heat source at $x=0$) from measuring the temperature,
$h(t)=w(l,t)$, and heat flux $f(t)=\rd_x w(x,t)|_{x=l} $ at point
$x=l$ for $0\le
 t \le T$. The solution to this problem, equation (\ref{goli}), can be derived
easily just by assuming the solution in the from
$w(x,t)=\sum^{\infty}_{j=0} a_{j}(t)(l-x)^j$ (e.g. see \cite{JRC}
chapter 2) as,

\be w(x,t)=\sum_{m=0}^{\infty} {\rd}^{\,m}_t h (t) \,{{(l-x)^{2
m}}\over(2 m)!}-\sum_{m=0}^{\infty} {\rd}^{\,m}_t f (t)\,{{(l-x)^{2
m+1}}\over(2 m+1)!}.\label{bertog}\ee

 Now we apply the results of Theorem (\ref{lmf4}) and Lemma (\ref{lemf6})
 to equation (\ref{babol}) to
 prove the following theorem.

\begin{thm}\label{lemf88} For $g(t)\in M[0,T]$ from equation (\ref{babol}) we have,\be
\rd_{x}\,u(x,t)|_{x=x_{0}}=-\rd^{1\over2}_{t}g(t).\label{trmp}\ee

\end{thm}
{\bf Proof.}  By applying $\rd^{1\over2}_{t}$ to both sides of
equation (\ref{babol}) we have, \be \rd^{1\over2}_{t} \int^{t}_{0}
\,K(x_0-x,t-\tau) \,u(x,\tau)\,d\tau =
\rd^{1\over2}_{t}g(t).\nn\hspace{1.3cm}(\ast)\ee

Assuming $u(x,t)$ as a function of $t$ is smooth on $[0,T]$, then
for the left hand side of above equation using (\ref{hop}) we have,
\bea \rd^{1\over2}_{t} \int^{t}_{0} \,K(x_0-x,t-\tau)
\,u(x,\tau)\,d\tau =\int^{t}_{0}
\big(\rd^{\,{1\over 2}}_{t-\tau}K(x_0-x,t-\tau)\big)u(x,\tau)d\tau \nn\\
= \int^{t}_{0}\,\rd_{x}\,K(x_0-x,t-\tau) \,u(x,\tau)\,d\tau \nn\\ =
-\int^{t}_{0}\,K(x_0-x,t-\tau)
\,\big(\rd_{x}\,u(x,\tau)\big)\,d\tau, \nn\hspace{2cm}(\ast\ast)\eea

where in the second line we used (\ref{hop2}) and in the third line
we used, \be\int^{t}_{0}\,\rd_{x}\,K(x_0-x,t-\tau)
\,u(x,\tau)\,d\tau +\int^{t}_{0}\,K(x_0-x,t-\tau)
\,\big(\rd_{x}\,u(x,\tau)\big)\,d\tau=0,\nn \ee which comes from
differentiating (\ref{babol}) with respect to variable $x$. Thus
from ($\ast$) and ($\ast\ast$) we get,\be
\int^{t}_{0}\,K(x_0-x,t-\tau) \,\big(\rd_{x}\,u(x,\tau)\big)\,d\tau
=-\rd^{1\over2}_{t}g(t).\nn\ee

The above equations is valid for $x\in [0,x_0)$ so one can take the
limit $ x\to x^-_0$ and use equation (\ref{hop3}) to get, \bea
\lim_{x\to x^-_0}\int^{t}_{0}\,K(x_0-x,t-\tau)
\,\big(\rd_{x}\,u(x,\tau)\big)\,d\tau &=&\,-\rd^{1\over2}_{t}g(t)\nn\\
\rd_{x}\,u(x,t)|_{x=x_{0}}&=&\,-\rd^{1\over2}_{t}g(t).\nn\eea\qed

The mathematical consequence of Theorem \ref{lemf88} is that it
converts the problem of semi-infinite bar (\ref{shahla}) into,
 \bea
 \rd^2_x u(x,t) &=&\rd_t u(x,t), \hspace{1.2cm} 0\le x < x_0,\hspace{0.6cm} 0\le
 t \le T,
\nn \\ u(x_0,t) &=& g(t), \hspace{1.9cm}  0\le
 t \le T,\label{dastte}\\ \rd_x
u(x,t)|_{x=x_0} &=& - \rd_t^{1\over2}g(t) \hspace{1.4cm} 0\le
 t \le T, \nn\\ u(x,0) &=& 0, \hspace{2.4cm}
0\le x < x_0 .\nn\eea This is similar to problem of finite bar
(equation (\ref{goli}) for $w(x,t)$). So by doing the substitutions
\bea l &\to & x_{0}, \nn\\ h(t) &\to & g(t)\nn\\ f(t)  &\to &
-\rd^{1\over2}_{t}g(t)\nn\eea in expression for $w(x,t)$ (equation
 (\ref{bertog})), one can get the following answer for $u(x,t)$,
\bea u(x,t)&=&\sum_{n=0}^{\infty} {\rd}^{\,n}_t( g (t)) {{(x_0-x)^{2
n}}\over(2 n)!}+ \sum_{n=0}^{\infty} {\rd}^{\,n}_t( \rd_t^{1\over
2}g (t)) {{(x_0-x)^{2 n+1}}\over(2 n+1)!},\nn\\&=&
\sum_{n=0}^{\infty}\rd^{n\over 2}_{t}(g (t))
{{(x_0-x)^{n}}\over(n)!}.\label{marjan}\eea

Because of the  appearance of extra half derivative in
(\ref{marjan}) the existence Theorem in comparison with solution of
finite bar problem (\ref{bertog}) is a bit different. We discuss
this matter in the next part where we present a simple
existence Theorem for equation (\ref{shahla}).\\

\begin{Rema}\label{zari23}{\em Before moving to the next section here we just
notice that a similar method can be used to find analytical solution
of more general Cauchy problem, \bea
 \rd^2_x v(x,t) &=&\rd_t v(x,t), \hspace{1cm} 0\le x < \infty,\hspace{0.6cm} 0\le
 t \le T,
\nn \\ v(x_0,t) &=& g(t), \hspace{1.7cm} x_0 \in (0,\infty),
\hspace{0.5cm} 0\le
 t \le T,\label{afsoon}\\ v(x,0) &=& \psi(x), \hspace{1.6cm}
0\le x<\infty.\nn\eea

That is the inverse heat conduction of semi-infinite bar with
initial temperature distribution $v(x,0) = \psi(x)$. For $v(x,t)$,
the solution of this problem, it can be proved  (e.g. see \cite{JRC}
equation 4.1.1), \be \int^{t}_{0} \,K(x'-x,t-\tau)
\,v(x,\tau)\,d\tau+\int^{\infty}_{x} dy\,N(x',y,t)\,\psi(y)\,
=v(x',t),\hspace{1cm}x  < x'\nn,\ee

for $0 <t \le T$, where $N(x,y,t)=\Phi_1(x-y,t)-\Phi_1(x+y,t)$ and
\be \Phi_{1}({ x},t)=\begin{cases}{{\,\,1 \over \Large{({{4 \pi
 t}\, })^{1\over 2}}}}\,
\text{e}^{-{\,\,{ x}^2\over \Large{ {4  t}}}} & \mbox { if }\,\,\, {x}\in{\mathbb{R}}, \,\,\, t > 0, \\
0  & \mbox { if }\,\,\,  {x}\in{\mathbb{R}}, \,\,\, t = 0.
\end{cases}\nn\ee This leads to the following equivalent problem of Volterra
integral equation of first kind for $0 < x \le x_0 $ and  $0 < t\le
T$, \be \int^{t}_{0} \,K(x_0-x,t-\tau) \,v(x,\tau)\,d\tau
=g(t)-\int^{\infty}_{x} dy\,N(x_{0},y,t)\,\psi(y)\,,\hspace{1cm} 0
\le x < x_0,\label{fundu123}\ee

and the following answer for $x_0 < x$ and  $0 < t\le T$, \be
v(x,t)=\int^{t}_{0} \,K(x-x_0,\,t-\tau)
\,g(\tau)\,d\tau+\int^{\infty}_{x_{0}} dy\,N(x,y,t)\,\psi(y)\,
 ,\hspace{1cm}x_0 < x.\nn\ee

For Equation (\ref{fundu123}) we may apply half derivative and do
the same procedure as in Theorem (\ref{lemf88}) to get,
\bea-\int^{t}_{0}\,K(x_0-x,t-\tau)
\,\big(\rd_{x}\,v(x,\tau)\big)\,d\tau={\rd}^{1\over 2}_{t}\,\Big
(g(t)-\int^{\infty}_{x} dy\,N(x_{0},y,t)\,\psi(y)\Big),\nn\eea for
$0 \le x < x_0$ and $0<t\le T$. By taking the limit $ x\to x^-_0$
one find the heat flux at $x=x_0$ as, \be  \rd_x (v(x,t))|_{x=x_0}\,
=\, -{\rd}^{1\over 2}_{t}\Big (g(t)-\int^{\infty}_{x}
dy\,N(x_{0},y,t)\,\psi(y)\Big)|_{x=x_0},\hspace{0.75cm}0<t\le
T.\nn\ee By substituting $v(x_{0},t)$ and $\rd_{x}v(x_{0},t)$ with
$h(t)$ and $f(t)$ in (\ref{bertog}) respectively  one finds the
analytical solution for Cauchy problem (\ref{afsoon}).}
\end{Rema}

\subsection {Existence Theorem}\label{ans} A non-zero smooth
function $P(t)$ on $[0,T]$ which satisfies (\ref{moon}) can not be
analytic every where in $[0,T]$. Due to Cauchy estimate Theorem, the
useful condition, $ \exists M, R>0, \,\forall n\in \mathbb{Z}_{+},
\sup_{t\in[0,T]}{|Q^{(\,n)}(t)|}< M{{( n!)}\over R^n}$\,, is valid
and only valid, if $Q(t)$ is analytic in $[0,T]$ (where $
Q^{(\,n)}(t)=\rd^{n}_{t}Q(t)$). Thus one needs to extend the
definition of analyticity for smooth functions such that it includes
functions which satisfy (\ref{moon}) and also provide some useful
condition like above. This is classically done by the use of Geverey
classes (\cite{Gev}, \cite{Hol}, \cite{Hol2}). A special class of
interest for heat equation is Gevrey class two, $G^{2}([0,t])$,
called Holmgren functions, which can be defined by condition of,

\be \exists M, R>0, \,\forall n\in \mathbb{Z}_{+}, \hspace{1cm}
\sup_{t\in[0,T]}{|P^{(\,n)}(t)|}< M{{( 2 n)!}\over R^{\,2 n}}.
\nn\ee

Further more Holmgren functions can be classified according to
constants in the above condition. For the matter of consistency with
literature we bring definition of Holmgren classes from \cite{JRC}
(definition 2.2.1)
\begin{definition}\label{khalil}{{\bf \text{(Holmgren)}} For the positive constants ${\gamma_1},{\gamma_2}$ and $C$,
the Holmgren class $H({\gamma_1},{\gamma_2},C,t_0)$ is set of smooth
functions $\phi$ defined on $|t-t_0|<{\gamma_2}$ that satisfy \be
\,\forall n\in \mathbb{Z}_{+}, \hspace{1cm}\sup_{|t-t_0|<{\gamma_2}}
{|\phi^{(\,n)}(t)|}< C{{( 2 n)!}\over {\gamma_1}^{2 n}}. \nn\ee }
\end{definition}

In this paper we have chosen our interval such that $t$ start from
zero ($t\in [0,T]$), this is equivalent of choosing
$t_{0}={\gamma_2}$ in above definition. So for convenience we take
the notation, \be \mathcal
{H}_l({\gamma_1},{\gamma_2},C)\,:\,=\,H({\gamma_1},{\gamma_2},C,{\gamma_2}).
\ee Thus by definition $\phi(t)\in\mathcal
{H}_l({\gamma_1},{\gamma_2},C)$ is defined for
$t\in[0,2{\gamma_2}]$.

\begin{lem}\label{lemf603}
1-\,\,If $P(t)\in {H}_l({\gamma_1},{\gamma_2},C)$ then, \,$
\rd_{\,t} P(t)\in {H}_l\big(\sqrt{2\over
5}\,\,{\gamma_1},\,\,{\gamma_2}\,\,,C_1\big)$ where $C_1=\,{24\over
5}{\gamma_1}^{-2} C$\\
2-\,\,If  $R(t)\in M[0,T]$ and $R(t)\in
{H}_l({\gamma_1},{\gamma_2},C)$ then, \,$ \rd^{1\over2}_{\,t}
R(t)\in {H}_l\big(\sqrt{2\over
5}\,\,{\gamma_1},\,\,{\gamma_2}\,\,,C_2\big)$ where $C_2={48
\sqrt{2\gamma_2}\over {5 \sqrt{\pi}}} {\gamma_1}^{-2} C$
\end{lem}

{\bf Proof.} Let us first show the following
 inequality for $1 \le \beta
\le {12\over 5}$ , \be (2j+1)(2j+2)\le 2 \beta\, {({{6}\over
\beta})^j},\hspace{0.5cm} j=0,1,2,...
 \nn\ee Assuming $(2j+1)(2j+2)\le c (d)^{j}$ then
the parametrization, $c=2\beta$, $d=6/\beta$ for $1 \le \beta
\le{12\over 5}$ can be derived matching the inequality for
$j=0,1,2$. The rest can be proved by induction.

Now for part 1 assuming $P(t)\in {H}_l({\gamma_1},{\gamma_2},C)$
then for $Q(t)=\rd_{t} P(t)$ one gets, \be
\sup_{t\in[0,T]}|Q^{(n)}(t)|=\sup_{t\in[0,T]}|P^{(n+1)}(t)|\le C{{(
2 n+2)!}\over {\gamma_1}^{2 n+2}}= (C {\gamma_1}^{-2 }) {{( 2
n)!}\over {\gamma_1}^{2 n}}(2n+1)(2n+2).\nn\ee Now by using the
above inequality, for $1 \le \beta \le {12\over 5}$ we have, \be
\sup_{t\in[0,T]}|Q^{(n)}(t)|\le (C {\gamma_1}^{-2 }) {{( 2 n)!}\over
{\gamma_1}^{2 n}}(2\beta)(6/\beta)^n= 2 \beta\, {\gamma_1}^{-2} C
{{( 2 n)!}\over {\gamma}^{2 n}},\nn\ee

where ${\gamma}={\gamma_1}\sqrt{\beta\over 6}$. Therefor
$\rd_{t}P(t)\in {H}_l\big(({\gamma_1}\sqrt{\beta\over
6}),\,{\gamma_2}\, ,(2 \beta\, {\gamma_1}^{-2} C)\big)$ for \,\,$1
\le \beta \le {12\over 5}$ and thus particularly for $\beta =
{12\over 5}$ we have, $ \rd_{t}P(t)\in {H}_l\big(\sqrt{2\over
5}\,\,{\gamma_1},\,\,{\gamma_2}\,\,,\,\,{24\over 5}{\gamma_1}^{-2} C
\big).$\\

For part 2, since $R(t)\in M[0,T]$ by using Theorem \ref{lmf4} part
one we can write, \bea |{\rd}^{n}_t ( \rd_t^{1\over 2}R
(t))|=|{\rd}^{1\over 2}_t ( \rd_t^{n}R (t))|&=&
{1\over\sqrt{\pi}}\,\,|\int^{t}_{0} d\tau\,{( \rd_t^{n+1}R
(\tau))\over \sqrt{t-\tau}}|\nn\\&=&|\int^{t}_{0}
d\tau\,{2{\sqrt{t-\tau}}\over
\sqrt\pi}(\rd_{\tau}^{\,n+2}\,R(\tau))|. \nn\eea

But since $\sqrt{t-\tau}\le \sqrt{2 \gamma_2}$\,for $t\in[0,2
\gamma_2]$ therefore,

\be |{\rd}^{n}_t ( \rd_t^{1\over 2}R (t))|\le|\int^{t}_{0}
d\tau\,{{2\sqrt{2\gamma_2}}\over \sqrt\pi}(\rd_{\tau}^{n+2} R
(\tau))|={2\sqrt{2\gamma_2}\over \sqrt\pi}|\rd_{t}^{n+1}
R(t)|.\nn\ee   Since $R(t)\in {H}_l({\gamma_1},{\gamma_2},C)$ by
part 1 we have, $\rd_t \, R(t) \in {H}_l\big({\gamma_1}\sqrt{2\over
5}\,,\,\,{\gamma_2}\,\,,\,\,{24\over 5}{\gamma_1}^{-2} C\big)$ thus,
$|\rd^{\,n+1}_t R(t)|< {24\over 5} \,C{\gamma_1}^{-2}({{( 2
n)!}\over {\gamma'_1}^{2 n}})$ where $\gamma'_1=\gamma_1
\sqrt{2\over 5}\,$. By inserting this in to the last line we get,

\be |{\rd}^{n}_t ( \rd_t^{1\over 2}g (t))|< \,\big({48
\sqrt{2\gamma_2}\over {5 \sqrt{\pi}}} {\gamma_1}^{-2} \,C
\,\big){{(2 n)}!\over{\gamma'_1}^{2 n}}.\nn\ee\qed

\begin{thm}\label{lemf621} If $g(t)\in
{H}_l({\gamma_1},\gamma_2,C)$ and $g(t)\in M[0, 2 \gamma_2]$ then
the power series, \be u(x,t)=\sum_{n=0}^{\infty} {\rd}^{\,n}_t( g
(t)) {{(x_0-x)^{2 n}}\over(2 n)!}+ \sum_{n=0}^{\infty} (
\rd_t^{n+{1\over 2}}g (t)) {{(x_0-x)^{2 n+1}}\over(2 n+1)!}, \nn\ee
converges uniformly and absolutely for $|x_0-x|\le r <\sqrt{2\over
5}\,{\gamma_1}$ and $u(x,t)$ is the solution of equation
(\ref{shahla}) for $t\in[0,2\gamma_2]$ and $|x_0-x|\le r <
\sqrt{2\over 5}\,{\gamma_1}$.
\end{thm}
{\bf Proof.} Since $g(t)$ belongs to ${H}_l({\gamma_1},\gamma_2,C)$,
it follows from Lemma (\ref{lemf603}) that $\rd^{1\over2}_{\,t}
g(t)\in {H}_l\big(\sqrt{2\over
5}\,\,{\gamma_1},\,\,{\gamma_2}\,\,,C_2\big)$. Therefore both $g(t)$
and $\,\rd^{1\over2}_{\,t} g(t)$ are in ${H}_l\big(\sqrt{2\over
5}\,\,{\gamma_1},\,\,{\gamma_2}\,\,,C_3\big)$ where $C_3=\max\{C,
C_2\}$. Considering equation ($\ref{hope}$), the series for $u(x,t)$
can be rewritten as, \be u(x,t)=\sum_{n=0}^{\infty}
\big({\rd}^{\,n}_t( g (t)) {{(x_0-x)^{2 n}}\over(2 n)!}+{\rd}^{n}_t
( \rd_t^{1\over 2} g (t)) {{(x_0-x)^{2 n+1}}\over(2
n+1)!}\big).\nn\ee The rest of proof is similar to
Holmgren result for problem (\ref{goli}) (e.g. see \cite{JRC} Theorem (2.3.1)). \qed\\

It should be noted that a solution $u(x,t)$ with the conditions in
Theorem \ref{lemf621} covers the whole $x\in[0,x_0]$ when, \be x_{0}
< {\sqrt{2\over 5}\,{\gamma_1}},  \ee and thus one can retrieve,
$u(0,t)$, the unknown heat source at $x=0$. Otherwise,

\be x_{0} \ge  {\sqrt{2\over 5}\,{\gamma_1}} \Longrightarrow u(x,t)
{\text{ is convergent for,\,\,\,\,\,}} x_0-{\sqrt{2\over
5}\,{\gamma_1}}< \,\, x \,\,\le \,\, x_0.\nn\ee\\

There are also functions like,
\be \psi(t)=\begin{cases}\,\text{exp}({{-1}\over t^{\,2}}) & \mbox { if }\,\,\, 0 < t \le T, \nn\\
0  & \mbox { if }\,\,\,  t = 0,
\end{cases}\ee
which satisfy extra property that for each $\gamma_1>0$ there exists
a $C_1=C_1(\gamma_1)>0$ such that $\psi\in
{H}_l({\gamma_1},\gamma_2,C_1(\gamma_1))$ for fixed $\gamma_2$ (e.g.
see \cite{Wider} or \cite{JRC} section 2.4). Therefore for such
functions by Theorem \ref{lemf621}, the series,
$\sum_{n=0}^{\infty}\rd^{n\over 2}_{t}(\psi (t))
{{(x_0-x)^{n}}\over(n)!},$ always converges uniformly and absolutely
for any $x\in[0,x_0]$ and $t\in[0,2\gamma_2]$.\\\\

\subsection{An example via 1d direct heat equation}\label{example} {The relation
between heat flux, $u_x(x,t)$ and time semi-derivative of
temperature $u(x,t)$ has been noted practically in electrochemistry
(in the context of diffusion equation where $u(x,t)$ represents
density) and is justified using direct heat equation with special
boundary conditions (e.g. see \cite{olhm2}, \cite{pold}).\\
In the following we bring, a case of direct one dimensional heat
conduction equation of {\em infinite bar} ($x\in(-\infty,+\infty)$)
with an initial distribution, $g(x)$, being confined to $x<0$. In
other word the initial distribution is zero for $x\ge 0$, thus it
provides an example for testing formulas (\ref{trmp}) and
(\ref{marjan}), since the left part of infinite bar ($x<0)$ can be
thought as external source at $x=0$ for the right part ($x\ge 0$).
Our method here again is based on rigorous properties of heat kernel
with respect to fractional derivatives. For the record, in the
following we keep thermal diffusivity constant, $\kappa$, in the
heat conduction equation.

\begin{lem}\label{lemf798}
Considering the following initial value problem,
\be \begin{cases}{\rd_{t}u({{x}},t)-\kappa\rd^{\,2}_{{x}}\,u({{x}},t)=0}, & \mbox \,\,\, ({x}\in{\mathbb{R}}, \,\,\, t > 0), \\
u({{x}},0) = g({{x}}),  & \mbox \,\,\, ({ x}\in{\mathbb{R}}),
\end{cases}\label{hagh}\ee

where $g(x)$ is bounded continuous function on $\mathbb{R}$  ($g({
x})\in C^{b}(\mathbb{R})$) and \be g({x})=0,\hspace{1.2cm}\,\,\,\,
x\ge 0.\label{mehdi}\ee

then $u({ x},t)$ for $x > 0$,  as a function of $t$, is in
$M[0,+\infty]$ and we have,

\be \rd_x u({ x},t) ={-1\over \sqrt{\kappa}}\,\,\,
\rd_{t}^{1\over2}u(x,t)),\hspace{1.2cm}t>0,\,\,\,\, x >
0,\label{lov}\ee
 \be u(x,t) = \sum_{n=0}^{\infty}
{{(a-x)^{n}}\over(n)!} {\rd^{n\over 2}_{t}(u(a,t))\over
\kappa^{n\over 2}} \hspace{1.2cm}t>0,\,\,\,\, a
> 0,\,\, x\in(0,a],\label{marjan3}\ee
\be \int^{t}_{0} \,K(x'-x,t-\tau) \,u(x,\tau)\,d\tau
=u(x',t),\hspace{1cm}x'  > x > 0,\label{ghodi}\ee where $K(x,t)
={{\,\,x\over \Large{ {2\sqrt{\pi\, \kappa}\,\,\,
 t^{3\over2}}}}}\,\text{e}^{-{\,\,x^2\over \Large{ {4\,\kappa\,t}}}}.$
\end{lem}

{\bf Proof.} With condition $g({ x})\in C^{b}(\mathbb{R})$ the
solution for equation (\ref{hagh}) exists in
$C^{\infty}({\mathbb{R}}\times(0,\infty))$, given by, \be u({
x},t)=\int_{\mathbb{R}}\Phi_{1}({ x}-{ y}, t)\, g({y})\, d{y},
\nn\ee where,\be \Phi_1({ x},t)=\begin{cases}{{\,\,1 \over
\Large{({{4 \pi \kappa t}\, })^{1\over 2}}}}\,
\text{e}^{-{\,\,|{ x}|^2\over \Large{ {4 \kappa t}}}} & \mbox { if }\,\,\, { x}\in{\mathbb{R}}, \,\,\, t > 0, \\
0  & \mbox { if }\,\,\,  { x}\in{\mathbb{R}}, \,\,\, t = 0,
\end{cases}\nn\ee
(e.g. see \cite{Evans}, section 2.3, Theorem 1). It is easy to show
that the function $\Phi_1({x},t)$ for $x\ne0$ as a function of $t$
is in $M[0,+\infty]$ and we have,

\be \rd_{x}\Phi_{1}({ x}, t)={-1\over
\sqrt{\kappa}}\,\,\rd^{1\over2}_{\,t} \Phi_{1}({ x},
t),\hspace{0.75cm}x>0,\, t\ge 0,\hspace{2cm} (\star)\nn\ee

It should be noted that function $\Phi_1({x},t)$ does not belong to
$M[0,+\infty]$ for $x=0$ and also for $x<0$ the relation ($\star$)
changes sign. Considering the condition (\ref{mehdi}) for $g({x})$
we have, \be u({ x},t)=\int_{{ y}< 0}\Phi_{1}({ x}-{ y}, t)\, g({
y})\, d{ y}. \hspace{3.5cm} (\star\star)\nn\ee

Now if ${{x}}>0$ then, $\Phi_{1}({ x}-{y}, t)$ appearing in the
integral ($\star\star$), belongs to $M[0,+\infty]$ since $({ x}-{ y})>0$. It follows that for any integer $n\ge0$,\\
 \be \rd^n_{t} u(x,t)|_{t=0}=\int_{{ y}< 0} \rd^n_{t}
\Phi_{1}({ x}-{ y}, t)|_{t=0}\, g({ y})\, d{ y}=0,\hspace{1cm} x>0
.\nn\ee

 Thus for $x > 0$, $u({ x},t)$ as a function of $t$ is in
$M[0,+\infty]$.

For $x>0$ we have, \bea \rd_x u({ x},t) &=&
\,\int^{0}_{-\infty} \rd_x \Phi_{1}({ x}-{y}, t)\, g({ y})\, d{ y},\nn\\
&=&{-1\over \sqrt{\kappa}}\, \int^{0}_{-\infty}
\rd^{1\over2}_{t} \Phi_{1}({ x}-{y}, t)\, g({ y})\, d{ y},\nn\\
&=&{-1\over \sqrt{\kappa}}\,\int^{0}_{-\infty}\,dy \, {2\over
\Gamma({1\over2})}\big(\int^{t}_{0} d\tau
\,(t-\tau)^{1\over2}\rd^{\,2}_{\tau}\Phi_1(x-y,\tau)\big)g(y),
\nn\\&=&{-1\over \sqrt{\kappa}}\,{2\over
\Gamma({1\over2})}\int^{t}_{0} d\tau
(t-\tau)^{1\over2}\rd^{\,2}_{\tau}\big(\int^{0}_{-\infty}\,dy
\,\Phi_1(x-y,\tau)\big)g(y), \nn\\&=& {-1\over
\sqrt{\kappa}}\,{1\over \Gamma({1\over2})}\int^{t}_{0}
d\tau\,(t-\tau)^{-{1\over2}}\,\rd_{\tau} \big(\int^{0}_{-\infty}
\Phi_1(x-y,\tau)\big)g(y)\,dy,\nn
\\&=&{-1\over \sqrt{\kappa}}\,\,\,(\rd_t^{1\over 2}\,\, u(x,t)),\nn\eea where in
the first line we applied Leibniz integral rule, in the third line
we used equation (\ref{khaleghi}) for semi-fractional derivative of
$\Phi_{1}({ x}-{y}, t)$ and in the line four the integrations order
is
changed, using Fubini theorem.\\

For equation (\ref{marjan3}), for $t>0,\,\,\,\, a
> 0,\,\, x\in(0,a]$ we have,
\bea \sum_{n=0}^{\infty} {{(a-x)^{n}}\over(n)!} {\rd^{n\over
2}_{t}(u(a,t))\over \kappa^{n\over 2}}&=&\, \sum_{n=0}^{\infty}
{{(a-x)^{n}}\over(n)!} \Big(\int^{0}_{-\infty} \, d{ y}{\rd^{n\over
2}_{t}(\Phi_{1}({x}-{y}, t))\over \kappa^{n\over 2}} \, g({
y})\Big)|_{x=a},\nn.    \\&=&\, \int^{0}_{-\infty} dy
 \{\sum_{n=0}^{\infty} {{(x-a)^{n}}\over(n)!}\,\Big({\rd^{n}_{x}(\Phi_{1}({x}-{y}, t))}\Big)|_{x=a}\} \, g({
y}), \nn\\&=&\, \int^{0}_{-\infty} dy (\Phi_{1}({x}-{y}, t)) \, g({
y}) = u(x,t),\nn\eea where in the second line above we used the
relation ($\star$). The series appearing in the second line, is just
the Taylor series of function $\Phi_{1}({x}-{y}, t)$ around $x=a$
which is convergent (for $x\in(0,a]$) to itself since
$\Phi_{1}({x}-{y}, t)$ is an analytic function with respect to
 $x$.\\

 The last part, equation (\ref{ghodi}), is just brought to make sure that
 solution of 1d direct heat equation ($\star\star$), for right part of the line ($x>0$), satisfies the
 same integral relation of inverse heat equation
 of semi-infinite bar (\ref{shahla}). Considering solution ($\star\star$),
 it is enough to show,
\be \int^{t}_{0} \,K(x'-x,t-\tau) \,\Phi_{1}({ x}-{ y}, \tau)
\,d\tau =\,\Phi_{1}({ x'}-{ y}, t), \hspace{1cm} x'> x > y.
\hspace{1cm} (\star\star\star)\nn\ee Using $K(x,t)=-2 \kappa\,\rd_x
\Phi_{1}(x,t)$ and \be \nn \int^{t}_{t'} \,d\tau \,K(x'-x,t-\tau) \,
K({ x}-{ z}, \tau-t')=K(x'-x,t-t'),\ee by starting from left hand
side we have, for  $x'> x > y$,
 \bea &&\int^{t}_{0}\,K(x'-x,t-\tau)
\,\Phi_{1}({ x}-{ y}, \tau) \,d\tau \nn\\ &&={1\over
{2\kappa}}\int^{t}_{0} \,d\tau \,K(x'-x,t-\tau) \int^{y}_{-\infty}
dz\, K({ x}-{ z}, \tau) \nn \\ &&={1\over
{2\kappa}}\int^{y}_{-\infty}
dz \int^{t}_{0} \,d\tau \,K(x'-x,t-\tau) \, K({ x}-{ z}, \tau) \nn\\
&&={1\over {2\kappa}}\int^{y}_{-\infty} dz  K({ x'}-{ z}, t) =
\Phi_{1}({ x'}-{ y}, t).\nn\eea\qed

\section{Conclusion}
In this paper we considered solving Volterra integral equation of
first kind (\ref{babol}). We examined some properties of function
space $M[0,T]$ with regards to fractional integrals and derivatives
relevant to this problem. In Theorem \ref{lmf4} we found relation
(\ref{hop}) which is an extension of Leibnitz integral rule to
derivatives of fractional order in a special case. In Theorem
\ref{lemf88}, we used results of previous section to deduce a new
boundary condition (\ref{trmp}) for equivalent formulation of
problem in term of Cauchy problem of parabolic partial differential
equation (\ref{shahla}), by which the solution is found by the
series (\ref{marjan}). In Theorem \ref{lemf621}, an existence
theorem is found for solution (\ref{marjan}) using Holmgren classes.
In section \ref{example} we provided an example for this problem
through a special case of initial value problem given in Lemma
\ref{lemf798}.

\end{document}